\documentclass[12pt, a4paper, twoside]{article}
\usepackage{amssymb}
\usepackage[leqno]{amsmath}
\usepackage{latexsym}
\usepackage{enumerate}
\usepackage{amsmath, amssymb, amscd}
\Roman{equation}

\DeclareMathOperator{\Pic}{Pic}
\DeclareMathOperator{\lcm}{lcm}
\DeclareMathOperator{\Cliff}{Cliff}

\begin{document}

\title{\textbf{\Large{On the splitting of Lazarsfeld-Mukai bundles on K3 surfaces II}}}

\author{Kenta Watanabe \thanks{Nihon University, College of Science and Technology,   7-24-1 Narashinodai Funabashi city Chiba 274-8501 Japan , {\it E-mail address:goo314kenta@mail.goo.ne.jp}, Telephone numbers: 090-9777-1974} }

\date{}

\maketitle 


\begin{abstract}

\noindent In this paper, we say that a rank 2 bundle splits if it is given by an extension of two line bundles. In the previous works, we gave a necessary condition for Lazarsfeld-Mukai bundles of rank 2 to split, under a numerical condition ([W2], Theorem 3.1). We gave the splitting types of them on a smooth quartic hypersurface in $\mathbb{P}^3$ ([W2], Proposition 3.1) as a corollary of it. However, the assertion of it contains a few mistakes. In this paper, we correct them, and give an application of the results in [W2].

\end{abstract}

\section{Introduction} Let $X$ be a K3 surface, $C$ be a smooth curve of genus $g\geq3$, and $Z$ be a base point free divisor on $C$ with $h^0(\mathcal{O}_C(Z))=r\;(r\geq2)$. Then the Lazarsfeld-Mukai bundle $E_{C,Z}$ of rank $r$ associated with $C$ and $Z$ is defined as the dual of the kernel $F_{C,Z}$ of the evaluation map
$${\rm{ev}}:H^0(\mathcal{O}_C(Z))\otimes\mathcal{O}_X\rightarrow \mathcal{O}_C(Z).$$
By the construction of it, if $|K_C\otimes\mathcal{O}_C(-Z)|$ is not empty, $E_{C,Z}$ is globally generated off the base points of it. Moreover, we have $h^1(E_{C,Z})=h^2(E_{C,Z})=0$, by easy computation. Conversely, it is known that if a vector bundle $E$ of rank $r$ with $h^1(E)=h^2(E)=0$ is globally generated, then there exist a smooth curve $C$ and a base point free divisor $Z$ on $C$ with $\dim |Z|=r-1$ such that $E=E_{C,Z}$ (for example [C-P], Lemma 1.2). In particular, if a rank 2 vector bundle $E$ is given by the extension
$$0\longrightarrow M\longrightarrow E\longrightarrow N\longrightarrow 0$$
of two non-trivial and base point free line bundles $N$ and $M$ satisfying $h^1(N)=h^1(M)=0$, then it is of type $E_{C,Z}$. Therefore, the problem of when a Lazarsfeld-Mukai bundle $E_{C,Z}$ of rank 2 is given by an extension of two line bundles is natural and interesting. For simplicity, in this paper, we say that a rank 2 bundle splits if it is given by an extension of two line bundles. Donagi and Morrison proved that if $E_{C,Z}$ is not simple, then the following assertion holds.

\newtheorem{thm}{Theorem}[section]

\begin{thm} {\rm{([D-M], Lemma 4.4)}}. Let $X$, $C$, and $Z$ be as above. Assume that $h^0(\mathcal{O}_C(Z))=2$. If $E_{C,Z}$ is not a simple bundle, then there exist two line bundles $M$ and $N$ on $X$ and a subscheme $Z^{'}\subset X$ of finite length such that:

$\;$

{\rm{(a)}} $h^0(M)\geq2,\;h^0(N)\geq2;$ 

\smallskip

\smallskip

{\rm{(b)}} $N$ is base point free;

\smallskip

\smallskip

{\rm{(c)}} There exists an exact sequence
$$0\longrightarrow M\longrightarrow E_{C,Z}\longrightarrow N\otimes\mathcal{J}_{Z^{'}}\longrightarrow0.$$
Moreover, if $h^0(M\otimes N^{\vee})=0$, then the length of $Z^{'}$ is zero.

\end{thm}

\noindent The exact sequence as in Theorem 1.1 is called Donagi-Morrison's extension. Moreover, if $h^0(M\otimes N^{\vee})=0$, then we have $E_{C,Z}\cong N\oplus M$. However, Theorem 1.1 does not give an answer of the above problem, in the case where $E_{C,Z}$ is simple.

In the previous work, we remarked that if $C$ is a very ample smooth curve on $X$, the Lazarsfeld-Mukai bundle $E_{C,Z}$ of rank 2 defined by $C$ and a base point free pencil $|Z|$ on $C$ is ACM and initialized with respect to $H=\mathcal{O}_X(C)$ ([W2], Proposition 2.3). Moreover, we have characterized a necessary condition for $E_{C,Z}$ to split, by ACM line bundles with respect to $H$, under a numerical condition ([W2], Theorem 3.1). Moreover, in the case where $X$ is a quartic hypersurface in $\mathbb{P}^3$ and $C$ is a hyperplane section of $X$, we have investigated the splitting type of $E_{C,Z}$, by using a numerical characterization of ACM line bundles on $X$ ([W2], Proposition 3.1). However, the statement of it contains some mistakes, and we indicated no concrete application of them. Therefore, in this paper, we correct the mistakes and give some applications of the previous results as in [W2].

Our plan of this paper is as follows. In section 2, we recall some basic facts about line bundles on K3 surfaces. In section 3, we recall several important theorems about the Clifford index of polarized K3 surfaces proved by Green and Lazarsfeld, and prepare some notations. In section 4, we review the previous results about the splitting of Lazarsfeld-Mukai bundles of rank 2 as in [W2], and correct the mistakes of the assertion. In section 5, we recall the existence theorem for polarized K3 surfaces which are dealt with in this paper, and prepare some lemmas to prove our main theorem. In section 6, we prove our main theorem.

$\;$

\noindent {\bf{Notations and conventions}}. We work over the complex number field $\mathbb{C}$. A curve and a surface are smooth projective. For a curve or a surface $Y$, we denote by $K_Y$ the canonical line bundle of $Y$ and denote by $|D|$ the linear system defined by a divisor $D$ on $Y$. 

For a curve $C$, the Clifford index of a line bundle $A$ on $C$ is definded as
$$\Cliff(A):=\deg(A)-2(h^0(A)-1).$$
The Clifford index of a curve $C$ is defined to be the minimum value of the Clifford index of line bundles on $C$. We denote it by $\Cliff(C)$. It is well known that
$$0\leq\Cliff(C)\leq\lfloor\displaystyle\frac{g-1}{2}\rfloor,$$
and a curve which has the maximal Clifford index is general in the moduli space of curves of genus $g$.

For a surface $X$, $\Pic(X)$ denotes the Picard lattice of $X$. A regular surface $X$ (i.e., a surface $X$ with $h^1(\mathcal{O}_X)=0$) is called a K3 surface if the canonical bundle $K_X$ of it is trivial. For a vector bundle $E$, we denote by $E^{\vee}$ the dual of it, and if we fix a very ample line bundle $H$ as a polarization on $X$, then we write $E\otimes H^{\otimes l}=E(l)$. We will say that a vector bundle $E$ is initialized with respect to a given polarization $H$, if 
$$H^0(E)\neq0,\text{ and }H^0(E(-1))=0.$$

\section{Linear systems on K3 surfaces}

In this section, we recall some classical results about the properties of base point free line bundles on K3 surfaces.

$\;$

\newtheorem{prop}{Proposition}[section]

\begin{prop} {\rm{([SD], Proposition 2.7)}}. Let $L$ be a numerical effective line bundle on a K3 surface $X$. Then $|L|$ is not base point free if and only if there exists an elliptic curve $F$, a smooth rational curve $\Gamma$ and an integer $k\geq2$ such that $F.\Gamma=1$ and $L\cong\mathcal{O}_X(kF+\Gamma)$.
\end{prop}

\smallskip

\smallskip

\begin{prop} {\rm{([SD], Proposition 2.6)}}. Let $L$ be a line bundle on a K3 surface $X$ such that $|L|\neq\emptyset$. Assume that $|L|$ is base point free. Then one of the following cases occurs.

$\;$

{\rm{(i)}} $L^2>0$ and the general member of $|L|$ is a smooth irreducible curve of genus $\displaystyle\frac{L^2}{2}+1$.

\smallskip

\smallskip

{\rm{(ii)}} $L^2=0$ and $L\cong\mathcal{O}_X(kF)$, where $k\geq1$ is an integer and $F$ is a smooth curve of genus one. In this case, $h^1(L)=k-1$.

\end{prop}

\noindent We note that, a linear system $|C|$ given by an irreducible curve $C$ with $C^2>0$ is base point free. Hence, by Proposition 2.2, any line bundle $L$ on a K3 surface has no base point outside its fixed components. On the other hand, by the classification of hyperelliptic linear systems on K3 surfaces which was given by B. Saint-Donat ([SD]), we have the following assertion.

\begin{prop}{\rm{(cf. [M-M], and [SD], Theorem 5.2)}}. Let $L$ be a numerical effective line bundle with $L^2\geq4$ on a K3 surface $X$. Then $L$ is very ample if and only if the following conditions are satisfied.

\smallskip

\smallskip

{\rm{(i)}} There is no irreducible curve $F$ with $F^2=0$ and $F.L\leq 2$.

\smallskip

{\rm{(ii)}} There is no irreducible curve $B$ with $B^2=2$ and $L\cong\mathcal{O}_X(2B)$.

\smallskip

{\rm{(iii)}} There is no $(-2)$-curve $\Gamma$ with $\Gamma.L=0$.

\end{prop}

\section{Clifford index of polarized K3 surfaces}

In this section, we recall several works about the Clifford index of polarized K3 surfaces, and prepare some notations. First of all, we remark the following theorem.

\begin{thm} {\rm{([G-L])}}. Let $H$ be a base point free and big line bundle on a K3 surface $X$ of sectional genus $g$. Then, for all smooth irreducible curve $C\in |H|$, the Clifford index of $C$ is constant. Moreover, if $\Cliff(C)<\lfloor\displaystyle\frac{g-1}{2}\rfloor$, then there exists a line bundle $L$ on $X$ such that $L\otimes\mathcal{O}_C$ computes the Clifford index of $C$ for any smooth irreducible curve $C\in |H|$.\end{thm}

\noindent By Theorem 3.1, for any base point free and big line bundle $H$ on a K3 surface $X$, we can define the Clifford index of it by the Clifford index of the curves belonging to the linear system $|H|$. Hence, we denote it by $\Cliff(H)$. If a given polarized K3 surface $(X,H)$ is not general (i.e., $\Cliff(H)<\lfloor\displaystyle\frac{g-1}{2}\rfloor$), then we can choose the line bundle $L$ as in Theorem 3.1 so that it satisfies the properties as in the following theorem.

\begin{thm} {\rm{([Kn], Theorem 8.3)}}. Let $X$, $H$, and $g$ be as in Theorem 3.1, and assume that $\Cliff(H)=c$. If $c<\lfloor\displaystyle\frac{g-1}{2}\rfloor$, then there exists a smooth curve $D$ on $X$ satisfying $0\leq D^2\leq c+2,\;2D^2\leq D.L$, and 
$$c=\Cliff(\mathcal{O}_X(D)\otimes\mathcal{O}_C)=D.H-D^2-2,$$
for any smooth curve $C\in |H|$.\end{thm}

\noindent Let $X$ and $H$ be as in Theorem 3.2. Then we set 
$$\mathcal{A}(H):=\{L\in\Pic(X)\;|\;h^0(L)\geq 2\text{ and }h^0(H\otimes L^{\vee})\geq2\}.$$
If $\mathcal{A}(H)\neq\emptyset$, then we set the integer $\mu(H)$ as follows.
$$\mu(H):=\min\{L.(H\otimes L^{\vee})-2\;|\;L\in\mathcal{A}(H)\}.$$
Then we set
$$\mathcal{A}^0(H):=\{L\in\mathcal{A}(H)\;|\; L.(H\otimes L^{\vee})=\mu(H)+2\}.$$
Since the smooth curve $D$ as in Theorem 3.2 satisfies
$$h^0(\mathcal{O}_X(D)|_C)=h^0(\mathcal{O}_X(D))\text{ and }h^0(H\otimes\mathcal{O}_X(-D)|_C)=h^1(\mathcal{O}_X(D)|_C),$$
for any smooth curve $C\in |H|$, we have the following assertion.

\newtheorem{cor}{Corollary}[section]

\begin{cor}. Let $X$ and $H$ be as in Theorem 3.2. Then
$$\Cliff(H)=\min\{\mu(H),\;\lfloor\displaystyle\frac{g-1}{2}\rfloor\}.$$
\end{cor}

\noindent Moreover, it is known that a line bundle $L$ belonging to $\mathcal{A}^0(H)$ has the following properties.

\begin{prop} {\rm{([J-K], Proposition 2.6)}}. Let $(X,H)$ be a polarized K3 surface of sectional genus $g\geq2$ such that $\mathcal{A}(H)\neq\emptyset$. Then $\mu(H)\geq0$ and any line bundle $L\in\mathcal{A}^0(H)$ satisfies the following properties.

$\;$

{\rm{(i)}} The {\rm{(}possibly empty\rm{)}} base divisor $\Delta$ of $L$ satisfies $H.\Delta=0$.

\smallskip

{\rm{(ii)}} $h^1(L)=0.$ \end{prop}

\section{Splitting of Lazarsfeld-Mukai bundles} 

In this section, we review previous works in [W2]. Let $X$ be a K3 surface, and let $H=\mathcal{O}_X(1)$ be a given very ample line bundle on $X$. Then we say that a vector bundle $E$ is arithmetically Cohen-Macaulay (ACM for short) with respect to $H$ if $H^1(E(l))=0$ for any integer $l\in\mathbb{Z}$. We can easily see that if we let $C\in|H|$ be a smooth curve, then the Lazarsfeld-Mukai bundle $E_{C,Z}$ of rank 2 defined by a base point free divisor $Z$ on $C$ with $h^0(\mathcal{O}_C(Z))=2$ is ACM with respect to $H$ ([W2], Proposition 2.3). By using this fact, we can obtain the following necessary condition for $E_{C,Z}$ to split, under a cohomological condition.

\begin{prop} {\rm{([W2], Theorem 3.1)}}. Let $X$, $C$ and $Z$ be as above. Assume that $L$ is a line bundle on $X$ such that $h^1(H\otimes L^{\vee})=0$ and $E_{C,Z}$ fits the exact sequence
$$0\longrightarrow H\otimes L^{\vee}\longrightarrow E_{C,Z}\longrightarrow L\longrightarrow 0.$$
Then:

$\;$

{\rm{(i)}} If $(H\otimes L^{\vee})^2<0$, then $h^0(L\otimes\mathcal{O}_C(-Z))\neq0.$

\smallskip

\smallskip

{\rm{(ii)}} $L$ is initialized with respect to $H$, and $h^0(L)\geq2$.

\smallskip

\smallskip

{\rm{(iii)}} $H\otimes L^{\vee}$ and $L$ are ACM with respect to $H$.

\end{prop}

\noindent Conversely, by the following assertions, we have a sufficient condition for $E_{C,Z}$ to 
split.

\begin{prop}. Let the notations be as in Proposition 4.1. If a line bundle $L$ on $X$ satisfies $L\otimes\mathcal{O}_C\cong\mathcal{O}_C(Z)$ and $h^1(H\otimes L^{\vee})=0$, then:

\smallskip

\smallskip

\noindent {\rm{(i)}} $|L|$ is an elliptic pencil.

\smallskip

\noindent {\rm{(ii)}} There exists the following exact sequence.
$$0\longrightarrow H\otimes L^{\vee}\longrightarrow E_{C,Z}\longrightarrow L\longrightarrow 0$$
 \end{prop}
 
 \noindent {\bf{Proof}}. Let $g$ be the genus of $C$. Note that since $C$ is not hyperelliptic, $g\geq3$. (i) Since $|Z|$ is a pencil, by the Riemann-Roch theorem, we have
$$g+1-\deg Z=h^1(\mathcal{O}_C(Z))\geq0.$$
Since 
$$H.(L\otimes H^{\vee})=\deg Z-(2g-2)\leq 3-g\leq0,$$
by the ampleness of $H$, we have $h^0(L\otimes H^{\vee})=0$. Since $L|_C\cong\mathcal{O}_C(Z)$, we have the following exact sequence.
$$0\longrightarrow L\otimes H^{\vee}\longrightarrow L\longrightarrow \mathcal{O}_C(Z)\longrightarrow 0$$
Moreover, since $h^1(H\otimes L^{\vee})=0$, we have $h^0(L)=h^0(\mathcal{O}_C(Z))=2$. Since $\mathcal{O}_C(Z)$ is base point free and $H$ is ample, $L$ is base point free. Therefore, by Proposition 2.2, $|L|$ is an elliptic pencil.

\smallskip

\smallskip

\noindent (ii) Since $h^1(H\otimes L^{\vee})=0$ and the exact sequence
$$0\longrightarrow H^0(\mathcal{O}_C(Z))^{\vee}\otimes L\otimes H^{\vee}\longrightarrow L\otimes H^{\vee}\otimes E_{C,Z}\longrightarrow \mathcal{O}_C\longrightarrow 0,$$
we have
$$h^0(L\otimes H^{\vee}\otimes E_{C,Z})=h^0(\mathcal{O}_C)=1.$$
Hence, we have $H\otimes L^{\vee}\subset E_{C,Z}$. Let $M$ be a saturated subsheaf of $E_{C,Z}$ satisfying $H\otimes L^{\vee}\subset M$. 

Assume that $h^0(M)=0$. Since $H\otimes L^{\vee}\subset M$, we have $H\otimes M^{\vee}\subset L$. On the other hand, by the exact sequence
$$0\longrightarrow M\longrightarrow E_{C,Z}\longrightarrow E_{C,Z}/M\longrightarrow 0,$$
we have
$$h^0(H\otimes M^{\vee})\geq h^0(E_{C,Z})=g+3-\deg Z\geq2.$$
By the assertion of (i), we have $L=H\otimes M^{\vee}$. Since $M=H\otimes L^{\vee}$ and $L$ is an elliptic pencil, we have
$$M.(H\otimes M^{\vee})=(H\otimes L^{\vee}).L=H.L=\deg Z.$$
Hence, we have the assertion.

Assume that $h^0(M)>0$. By the exact sequence
$$0\longrightarrow H^0(\mathcal{O}_C(Z))^{\vee}\otimes M^{\vee}\longrightarrow E_{C,Z}\otimes M^{\vee}\longrightarrow M^{\vee}\otimes K_C(-Z)\longrightarrow 0,$$
we have 
$$h^0(M^{\vee}\otimes K_C(-Z))\geq h^0(E_{C,Z}\otimes M^{\vee})>0.$$
Hence, we have $\deg (M^{\vee}\otimes K_C(-Z))=2g-2-\deg Z-M.H\geq0$. On the other hand, since
$$M.H\geq (H\otimes L^{\vee}).H=2g-2-\deg Z,$$
we have $M.H=(H\otimes L^{\vee}).H.$ Therefore, we have $M=H\otimes L^{\vee}$. By the same reason as above, we have the assertion.$\hfill\square$

\newtheorem{lem}{Lemma}[section]

\begin{lem} {\rm{([LC], Proof of Lemma 3.2)}}. Let $X$ be a K3 surface, $C$ be a smooth curve on $X$, and $Z$ be a base point free divisor on $C$ with $h^0(\mathcal{O}_C(Z))=2$ such that $|K_C\otimes\mathcal{O}_C(-Z)|\neq\emptyset$. Assume that $Q$ is a torsion free sheaf of rank one on $X$. Then if there exists a surjective morphism $\varphi:E_{C,Z}\longrightarrow Q$, then $Q^{\vee\vee}$ is base point free and not trivial.
\end{lem}

\noindent By Lemma 4.1, we have the following proposition.

\begin{prop}. Let the notations be as in Proposition 4.1, and let $\deg Z=d$. If $\Cliff(H)=d-2$ and a saturated subsheaf $L$ of $E_{C,Z}$ satisfies $L.H\geq d$, then $L$ is ACM and initialized with respect to $H$, and there exists the following exact sequence.
$$0\longrightarrow L\longrightarrow E_{C,Z}\longrightarrow H\otimes L^{\vee}\longrightarrow 0$$
\end{prop}

\noindent {\bf{Proof}}. Since $L\subset E_{C,Z}$ is saturated, we have $(E_{C,Z}/L)^{\vee\vee}=H\otimes L^{\vee}$. Since $L.(H\otimes L^{\vee})\leq d$, by the assumption, we have $L^2\geq L.H-d\geq0$. Hence, we have $h^0(L)\geq2$. Since $$h^0(K_C(-Z))=h^1(\mathcal{O}_C(Z))=g+1-d\text{ and }d\leq\lfloor\displaystyle\frac{g+3}{2}\rfloor,$$
by Lemma 4.1, $h^0(H\otimes L^{\vee})\geq2$. Hence, by Corollary 3.1, we have
$$\Cliff(H)=d-2\leq\mu(H)\leq L.(H\otimes L^{\vee})-2.$$
Therefore, we have $L.(H\otimes L^{\vee})=d$. Then we have the following exact sequence.
$$0\longrightarrow L\longrightarrow E_{C,Z}\longrightarrow H\otimes L^{\vee}\longrightarrow 0$$
By Proposition 3.1, we have $h^1(L)=0$. Hence, by Proposition 4.1, $L$ is ACM and initialized with respect to $H$.$\hfill\square$

$\;$
 
\noindent If $X$ is a quartic hypersurface in $\mathbb{P}^3$ and $C\subset X$ is a hyperplane section of $X$, then we have a complete numerical characterization of ACM line bundles with respect to $H$. 

\begin{prop} {\rm{([W1], Theorem 1.1)}}. Let $X$ be a smooth quartic hypersurface in $\mathbb{P}^3$. Let $C$ be a smooth hyperplane section of $X$, let $H=\mathcal{O}_X(C)$, and let $D$ be a non-zero effective divisor on $X$. Then the following conditions are equivalent.

\smallskip

\smallskip

\noindent {\rm{(i)}} $\mathcal{O}_X(D)$ is an ACM and initialized line bundle with respect to $H$.

\smallskip

\noindent {\rm{(ii)}} One of the following cases occurs.

\smallskip

\smallskip

{\rm{(a)}} $D^2=-2$ and $1\leq C.D\leq 3$.

\smallskip

{\em{(b)}} $D^2=0$ and $3\leq C.D\leq 4$.

\smallskip

{\rm{(c)}} $D^2=2$ and $C.D=5$.

\smallskip

{\rm{(d)}} $D^2=4,\;C.D=6$ and $|D-C|=|2C-D|=\emptyset$. \end{prop}

\noindent In the case where $X$ and $C$ are of as in Proposition 4.4, we get more concrete descriptions about the splitting type of $E_{C,Z}$. Although we have stated it in [W2] already, the assertion is not correct. Therefore, we correct it as follows.

\begin{prop} {\rm{([W2], Proposition 3.1)}}. Let $X$, $C$ and $H$ be as in Proposition 4.4. Moreover, let $Z$ be a base point free divisor on $C$ such that $h^0(\mathcal{O}_C(Z))=2$. Then, if there exists a line bundle $L$ with $h^1(H\otimes L^{\vee})=0$ such that $E_{C,Z}$ fits the exact sequence
$$0\longrightarrow H\otimes L^{\vee}\longrightarrow E_{C,Z}\longrightarrow L\longrightarrow 0,$$
then we get $(L.H,L^2)=(5,2), (3,0)$ or $(4,0)$.
\end{prop}

\noindent {\bf{Proof}}. We note that, by Propsition 4.1 (ii), we have $h^0(L)\geq2$. Hence, $L$ is not trivial. Since $|Z|$ is a pencil, we have $h^1(\mathcal{O}_C(Z))=4-\deg Z$. Since $C$ is a trigonal curve, we have $\deg Z=3 $ or 4. Since $h^1(H\otimes L^{\vee})=0$, by Proposition 4.1 (iii), $L$ is ACM. By Proposition 4.1 (ii), $L$ is initialized. Hence, by Proposition 4.4, we have
$$(L.H,L^2)=(1,-2),\;(2,-2),\;(3,0),\; (4,0),\text{ or }(5,2).$$
Since $H$ is very ample, if $L.H\leq2$, the movable part of $|L|$ is empty and hence, $h^0(L)=1$. However, this is a contradiction.$\hfill\square$

$\;$

\noindent By Proposition 4.4, the line bundle $L$ satisfying $(L.H,L^2)=(3,0)$ or (4,0) as in Proposition 4.5 is ACM with respect to $H$ and hence, $h^1(H\otimes L^{\vee})=0$. Hence, if such a line bundle $L$ satisfies $L\otimes \mathcal{O}_C\cong \mathcal{O}_C(Z)$, by Proposition 4.2, there exists an exact sequence as in Proposition 4.5.

\section{Existence theorem}

In this section, we recall the following existence theorem for a certain polarized K3 surface of Picard number 2.

\begin{thm}{\rm{(cf. [J-K], Proposition 4.2 and Lemma 4.3)}}. Let $d$ and $g$ be integers with $g\geq3$ and $3\leq d\leq \lfloor\displaystyle\frac{g+3}{2}\rfloor $. Then there exists a K3 surface $X$ with $\Pic(X)=\mathbb{Z}H\oplus\mathbb{Z}F$ such that $H$ is a base point free line bundle with $H^2=2g-2$ and the general member of $|F|$ is an elliptic curve with $H.F=d$. Moreover, $\Cliff(H)=d-2$.
\end{thm}

\noindent Note that, in Theorem 5.1, since $H$ is base point free, the general member $C$ of $|H|$ is a smooth curve, by Proposition 2.2. Moreover, the gonality of $C$ is $d$ and it can be computed by the pencil on $C$ which is given by the restriction of the elliptic pencil $|F|$ on $X$ ([J-K], proof of Lemma 4.3 and Theorem 4.4). We prepare the following proposition.

\begin{prop}. Let the notations be as in Theorem 5.1. Then the following statements hold.

\smallskip

\smallskip

{\rm{(i)}} $X$ contains a $(-2)$-curve if and only if $d|g$. Moreover, in this case, if we set $g=md$, then $H\otimes F^{\vee\otimes m}$ is a unique $(-2)$-vector up to sign in $\Pic(X)$ and the member of $|H\otimes F^{\vee\otimes m}|$ is a $(-2)$-curve.

\smallskip

\smallskip

{\rm{(ii)}} Let $L$ be a line bundle which is given by an elliptic curve on $X$. Then;

\smallskip

{\rm{(a)}} If $d|g$, then $L=F$.

{\rm{(b)}} If $d\mid \hspace{-.67em}/g$, then $L=F$ or there exist integers $m$ and $n$ such that 
$$n>0,\;n(g-1)+md=0,\;n(g-1)=\lcm(g-1,d),$$
$\text{ and }L=H^{\otimes n}\otimes F^{\otimes m}$.

\smallskip

\smallskip

{\rm{(iii)}} $H$ is very ample.

\smallskip

\smallskip

{\rm{(iv)}} If $d\leq g-1$, then $H\otimes F^{\vee}$ is base point free.

\end{prop}

\noindent {\bf{Proof}}. (i) By easy computation, $H\otimes F^{\vee\otimes m}$ with $g=md$ is a unique $(-2)$-vector up to sign in $\Pic(X)$ and this case occurs precisely when $d|g$. In this case, since $H.(H\otimes F^{\vee \otimes m})=g-2>0$, we have $|H\otimes F^{\vee\otimes m}|\neq\emptyset.$ By the uniqueness of the $(-2)$-vector, the member of $|H\otimes F^{\vee\otimes m}|$ is a $(-2)$-curve.

\smallskip

\smallskip

(ii) Let $L=H^{\otimes n}\otimes F^{\otimes m}$, and assume that $L$ is a line bundle which is given by an elliptic curve on $X$. Then $L^2=2n^2(g-1)+2nmd=0$. If $n=0$, then, by the assumption, we have $m=1$. Assume that $n\neq0$. Then we have $n(g-1)+md=0$. Since $L.F=nd\geq0$, we have $n>0$, and since the general member of $|L|$ is irreducible, we have $n(g-1)=\lcm(g-1,d)$. Moreover, this case occurs only if $d\mid \hspace{-.67em}/g$. In fact, if $d|g$, we set $k=\displaystyle\frac{g}{d}$. Then, by the statement of (i), the member of $|H\otimes F^{\vee\otimes k}|$ is a $(-2)$-curve and $L.(H\otimes F^{\vee\otimes k})=-n<0$. Hence, $|L|$ is not base point free. This is a contradiction.

\smallskip

\smallskip

(iii) First of all, we set $m=\displaystyle\frac{g}{d}$. Since $H.(H\otimes F^{\vee\otimes m})=g-2>0$, by the assertion of (i), $H$ is ample. Moreover, if $L$ is a line bundle given by an elliptic curve on $X$, by the statement of (ii), we have $L.H\geq d\geq3$. Assume that there exists a base point free line bundle $B$ of sectional genus 2 with $H\cong B^{\otimes 2}$. Since $H^2=8$, we have $g=5$. In this case, by easy computation, we can find a unique 2-vector $H\otimes F^{\vee}$ up to sign in $\Pic(X)$ precisely when $H.F=3$. Since $H.(H\otimes F^{\vee})=5>0$, we have $B\cong H\otimes F^{\vee}$. However, we have the contradiction $H\cong F^{\otimes 2}$. Therefore, by Proposition 2.3, $H$ is very ample.

\smallskip

\smallskip

(iv) In the case where $d\leq g-1$, we have $(H\otimes F^{\vee})^2=2g-2-2d\geq0$. Since $H.(H\otimes F^{\vee})=2g-2-d\geq g-1>0$, we have $|H\otimes F^{\vee}|\neq\emptyset.$ If $d\mid \hspace{-.67em}/g$, by the assertion of (i), there is no $(-2)$-curve on $X$. Hence, $|H\otimes F^{\vee}|$ is base point free. If $d|g$, then for $m=\displaystyle\frac{g}{d}$, the member of $|H\otimes F^{\vee\otimes m}|$ is a unique $(-2)$-curve on $X$, and since $(H\otimes F^{\vee}).(H\otimes F^{\vee\otimes m})=g-2-d>0$ and $F.(H\otimes F^{\vee\otimes m})=d\geq3$, by Proposition 2.1, we have the assertion.$\hfill\square$

\section{Main results}

In this section, we give some applications of previous works as in section 4.

\begin{thm}. Let $d$ and $g$ be integers with $g\geq3$ and $3\leq d\leq \lfloor\displaystyle\frac{g+3}{2}\rfloor$. Then let $X$ be a K3 surface with $\Pic(X)=\mathbb{Z}H\oplus\mathbb{Z}F$ such that $H$ is a base point free line bundle with $H^2=2g-2$ and the general member of $|F|$ is an elliptic curve with $H.F=d$. Moreover, let $C\in |H|$ be a smooth curve and $Z$ be a divisor on $C$ which gives a gonality pencil $|Z|$ on $C$, {\rm{(}}that is, $\deg Z=d${\rm{)}}. Then $E_{C,Z}$ has the following properties.

$\;$

\noindent {\rm{(i)}} Assume that $d=g\;{\rm{(}}i.e., d=g=3{\rm{)}}$, or $g-1$, and let $L$ be a line bundle on $X$. Then the following conditions are equivalent.

\smallskip

\smallskip

{\rm{(a)}} There exists the following exact sequence.
$$0\longrightarrow H\otimes L^{\vee}\longrightarrow E_{C,Z}\longrightarrow L\longrightarrow 0$$

\smallskip

{\rm{(b)}} $|L|$ is an elliptic pencil on $X$, and $\mathcal{O}_C(Z)=L|_C$.

$\;$

\noindent {\rm{(ii)}} If $d=g$ {\rm{(}}i.e., $d=g=3${\rm{)}}, $E_{C,Z}$ is $H$-slope-stable.

$\;$

\noindent {\rm{(iii)}} If $d=g-1$, the following conditions are equivalent.

\smallskip

\smallskip

{\rm{(a)}} $E_{C,Z}$ is not $H$-slope-stable.

\smallskip

{\rm{(b)}} $E_{C,Z}$ is strictly $H$-slope-semistable.

\smallskip

{\rm{(c)}} $\mathcal{O}_C(Z)\cong F|_C$ or $H\otimes F^{\vee}|_C$.

$\;$

\noindent {\rm{(iv)}} If $d<g-1$, then the following conditions are equivalent.

\smallskip

\smallskip

{\rm{(a)}} $E_{C,Z}$ is not $H$-slope-stable.

\smallskip

{\rm{(b)}} $E_{C,Z}$ is not $H$-slope-semistable.

\smallskip

{\rm{(c)}} $\mathcal{O}_C(Z)\cong F|_C$.

\end{thm}

\noindent First of all, we prove the following lemma to prove Theorem 6.1.

\begin{lem}. Let $E_{C,Z}$ be as in Theorem 6.1, and let $L$ be a line bundle such that $E_{C,Z}$ fits the following exact sequence.
$$0\longrightarrow H\otimes L^{\vee}\longrightarrow E_{C,Z}\longrightarrow L\longrightarrow 0$$
Then $L$ is a non-trivial and base point free line bundle satisfying 
$$h^1(L)=h^0(L\otimes H^{\vee})=0.$$
\end{lem}

\noindent {\bf{Proof}}. By the Riemann-Roch theorem, we have $$h^0(K_C\otimes\mathcal{O}_C(-Z))=g-d+1\geq1,$$
and hence, we have $|K_C\otimes\mathcal{O}_C(-Z)|\neq\emptyset$. Therefore, $E_{C,Z}$ is globally generated off the base points of $|K_C\otimes\mathcal{O}_C(-Z)|$. Assume that $E_{C,Z}$ fits the following exact sequence.
$$0\longrightarrow H\otimes L^{\vee}\longrightarrow E_{C,Z}\longrightarrow L\longrightarrow 0$$
Since $h^1(E_{C,Z})=h^2(E_{C,Z})=0$, if we assume that $L$ is not initialized with respect to $H$, then we have 
$$h^1(L)=h^2(H\otimes L^{\vee})=h^0(L\otimes H^{\vee})\neq0.$$
By Lemma 4.1, $L$ is globally generated and not trivial. Hence, by Proposition 2.2, there exist an elliptic curve $D$ on $X$ and an integer $r\geq2$ such that $L\cong\mathcal{O}_X(rD)$. However, by Proposition 5.1 (ii), we have the contradiction
$$d=L.(H\otimes L^{\vee})=L.H\geq rd>d.$$
Therefore, we have $h^1(L)=h^0(L\otimes H^{\vee})=0$.$\hfill\square$

$\;$

\noindent {\bf{Proof of Theorem 6.1}}. (i) (a)$\Rightarrow$(b). Assume that there exists a line bundle $L$ such that $E_{C,Z}$ fits the following exact sequence.
$$0\longrightarrow H\otimes L^{\vee}\longrightarrow E_{C,Z}\longrightarrow L\longrightarrow 0$$
We consider the case where $d=g=3$. If we let $\Gamma=H\otimes F^{\vee}$, then by Proposition 5.1 (i), the member of $|\Gamma|$ is a $(-2)$-curve. Since $\Pic(X)$ is generated by $H$ and $\Gamma$, we set $L\cong H^{\otimes n}\otimes\Gamma^{\otimes m}\;(m,n\in\mathbb{Z})$. Note that $m<0$. In fact, since $L$ is globally generated and not trivial, we have $h^0(L)\geq2$. We have $n\geq1$. Therefore, if $m\geq0$, then $L$ is not initialized. This contradicts Lemma 6.1. Since $H\otimes\Gamma^{\vee}=F$ has no base point and $L\cong H^{\otimes(n+m)}\otimes F^{\vee m}$, we have $L.F=3(n+m)\geq0$. Since $h^0(L\otimes H^{\vee})=0$, we must have $n+m=0$ and $h^1(L)=h^2(H\otimes L^{\vee})=0$. We have $m=-1$. Hence, $H\otimes L^{\vee}=\Gamma$. Since $(H\otimes L^{\vee})^2=-2<0$, by Proposition 4.1 (i), we have $h^0(L\otimes\mathcal{O}_C(-Z))>0$ and
$$\deg (L\otimes\mathcal{O}_C(-Z))=L.H-\deg Z=L.(H\otimes L^{\vee})-\deg Z=0.$$
Hence, we have $L|_C\cong \mathcal{O}_C(Z)$.

\smallskip

\smallskip

We consider the case where $d=g-1$. Let $L=H^{\otimes n}\otimes F^{\otimes m}$. Since $L.F\geq0$, we have $n\geq0$. Moreover, since $L$ is globally generated, we have
$$(2g-2)n(m+n)=L^2\geq0.$$
Since $h^1(L)=0$, if $n=0$, then we have $m=1$ and hence, $L\cong F$. If $n\geq1$, we have $m\geq-n$. Since $h^1(L)=0$, if $m=-n$, then we have $L\cong H\otimes F^{\vee}$. If $m>-n$, then there exists an injection $(H\otimes F^{\vee})^{\otimes (n-1)}\hookrightarrow L\otimes H^{\vee}$. Since $h^0(H\otimes F^{\vee})>0$, we have $h^0(L\otimes H^{\vee})\neq0$. Hence, by Lemma 6.1, this case does not occur. Therefore, by Proposition 5.1, $|H\otimes L^{\vee}|$ is an elliptic pencil. Hence, by the exact sequence
$$0\longrightarrow H^0(\mathcal{O}_C(Z))^{\vee}\otimes L\otimes H^{\vee}\longrightarrow E_{C,Z}\otimes L\otimes H^{\vee}\longrightarrow L\otimes\mathcal{O}_C(-Z)\longrightarrow 0,$$
we have $h^0(E_{C,Z}\otimes L\otimes H^{\vee})=h^0(L\otimes\mathcal{O}_C(-Z)).$ Since $H\otimes L^{\vee}\subset E_{C,Z}$, we have $h^0(L\otimes \mathcal{O}_C(-Z))>0$. Since
$$\deg(L\otimes\mathcal{O}_C(-Z))=L.H-\deg Z=L.(H\otimes L^{\vee})-\deg Z=0,$$
we have $L|_C\cong\mathcal{O}_C(Z)$.

\smallskip

\smallskip

\noindent (b)$\Rightarrow$(a). By Proposition 5.1 (ii), if $L$ is an elliptic pencil, then $L=F$ or $H\otimes F^{\vee}$. Since $h^1(H\otimes L^{\vee})=0$, by Proposition 4.2, we have the assertion.

$\;$

\noindent Before the proof of the latter part of Theorem 6.1, we prepare the following lemma.

\begin{lem}. Let the notations be as in Theorem 6.1. Let $L$ be a non-trivial line bundle on $X$. Then $L$ is ACM and initialized with respect to $H$ if and only if $L=F$ or $H\otimes F^{\vee}$.\end{lem}

{\bf{Proof}}. We assume that $L$ is ACM and initialized with respect to $H$. If we set $L=H^{\otimes n}\otimes F^{\otimes m}$, then $h^1(L\otimes H^{\vee\otimes n})=0$. Hence, we have $m=-1,0$ or 1. Since $L$ is non-trivial, we have $L=F$ or $H\otimes F^{\vee}$. Conversely, we assume that $L=F$ or $H\otimes F^{\vee}$. By Proposition 5.1, $F$ is initialized with respect to $H$. Moreover, the general member of $|H\otimes F^{\vee}|$ is a smooth curve. Hence, for any $l\in\mathbb{Z}$ such that $l\geq0$, we have $h^1(H^l\otimes F^{\vee})=0$. Obviously, since $h^1(H^l\otimes F)=0$ for any $l\geq0$, $F$ is ACM with respect to $H$. Hence, $H\otimes F^{\vee}$ is also ACM. Since $h^0(H\otimes F^{\vee})>0$, it is initialized with respect to $H$.$\hfill\square$

$\;$

\noindent Now we return to the proof of Theorem 6.1.

\smallskip

\smallskip

(ii) We consider the case where $d=g=3$. We assume that there exists a saturated sub-line bundle $L$ of $E_{C,Z}$ such that $L.H\geq\displaystyle\frac{H^2}{2}=2$. By Lemma 4.1, $H\otimes L^{\vee}=(E_{C,Z}/L)^{\vee\vee}$ is base point free. Hence, by Proposition 2.3 and the Hodge index theorem, we have $H.(H\otimes L^{\vee})\geq3$. However, this means that $L.H\leq1$. This is a contradiction.

$\;$

(iii) We consider the case where $d=g-1$. (b)$\Rightarrow$(a) is clear. 

\smallskip

\smallskip

\noindent (a)$\Rightarrow$(c). Assume that $E_{C,Z}$ is not $H$-slope-stable. Let $M\subset E_{C,Z}$ be a saturated sub-line bundle with $M.H\geq\displaystyle\frac{H^2}{2}=g-1$. Since $d=g-1$, by Proposition 4.3, $M$ is ACM and initialized with respect to $H$, and there exists the following exact sequence.
$$0\longrightarrow M\longrightarrow E_{C,Z}\longrightarrow H\otimes M^{\vee}\longrightarrow 0$$
Since $h^1(M)=0$, by the exact sequence
$$0\longrightarrow H^0(\mathcal{O}_C(Z))^{\vee}\otimes M^{\vee}\longrightarrow M^{\vee}\otimes E_{C,Z}\longrightarrow M^{\vee}\otimes K_C(-Z)\longrightarrow 0,$$
we have 
$$h^0(M^{\vee}\otimes K_C(-Z))=h^0(M^{\vee}\otimes E_{C,Z})>0.$$
By Lemma 6.2, we have $M=F$ or $H\otimes F^{\vee}$. Therefore, we have
$$\deg(M^{\vee}\otimes K_C(-Z))=2g-2-d-M.H=0.$$
Since $H\otimes M^{\vee}=H\otimes F^{\vee}$ or $F$, we have $\mathcal{O}_C(Z)=H\otimes F^{\vee}|_C$ or $F|_C$.

\smallskip

\smallskip

\noindent (c)$\Rightarrow$(b). Assume that $\mathcal{O}_C(Z)=F|_C$ or $H\otimes F^{\vee}|_C$. Let $L=F$ or $H\otimes F^{\vee}$. Since $h^1(H\otimes L^{\vee})=0$, by Proposition 4.2, we have the following exact sequence.
$$0\longrightarrow H\otimes L^{\vee}\longrightarrow E_{C,Z}\longrightarrow L\longrightarrow 0$$
Assume that $E_{C,Z}$ is not $H$-slope semi-stable. Since the maximal destabilizing sheaf $M$ of $E_{C,Z}$ satisfies $M.H>g-1=d$, by Proposition 4.3, $M$ is initialized and ACM with respect to $H$. Therefore, by Lemma 6.2, we have $M=F$ or $H\otimes F^{\vee}$. However, since $M.H=d$, this is a contradiction. Hence, $E_{C,Z}$ is strictly $H$-slope semi-stable.

$\;$

(iv) We consider the case where $d<g-1$. (b)$\Rightarrow$(a) is clear.

\smallskip

\smallskip

\noindent (a)$\Rightarrow$(c). We assume that there exists a saturated sub-line bundle $M\subset E_{C,Z}$ such that $M.H\geq\displaystyle\frac{H^2}{2}=g-1$. By the same way of the proof of (iii) (a)$\Rightarrow$(c), we have $M=F$ or $H\otimes F^{\vee}$. Since $d<g-1$, we have $M=H\otimes F^{\vee}$. Since $h^1(H\otimes F^{\vee})=0$, by the exact sequence
$$0\longrightarrow H^0(\mathcal{O}_C(Z))^{\vee}\otimes F\otimes H^{\vee}\longrightarrow F\otimes H^{\vee}\otimes E_{C,Z}\longrightarrow F\otimes H^{\vee}\otimes K_C(-Z)\longrightarrow 0,$$
we have
$$h^0(F\otimes H^{\vee}\otimes K_C(-Z))=h^0(F\otimes H^{\vee}\otimes E_{C,Z})>0.$$
On the other hand, since $\deg(F\otimes H^{\vee}\otimes K_C(-Z))=0$, we have $\mathcal{O}_C(Z)\cong F|_C$.

\smallskip

\smallskip

\noindent (c)$\Rightarrow$(b). We assume that $\mathcal{O}_C(Z)=F|_C$. Since $h^1(H\otimes F^{\vee})=0$, by Proposition 4.2, we have $H\otimes F^{\vee}\subset E_{C,Z}$. Since $d<g-1$, we have $H.(H\otimes F^{\vee})=2g-2-d>g-1$. Hence, $E_{C,Z}$ is not $H$-slope semi-stable.$\hfill\square$

$\;$

\noindent In Theorem 6.1, the case where $d=g-1$ occurs precisely when $(g,d)=(4,3)$ or (5,4). In particular, if $(g,d)=(4,3)$, then we have the following assertion.

\begin{cor}. Let the notations be as in Theorem 6.1. If $(g,d)=(4,3)$, then $W^1_3(C)=\{F|_C,\;H\otimes F^{\vee}|_C\}$.\end{cor}

\noindent {\bf{Proof}}. Let $Z$ be a divisor on $C$ with $h^0(\mathcal{O}_C(Z))=2$ and $\deg(Z)=3$. Since $h^0(E_{C,Z})=4$, we have $\dim\bigwedge^2H^0(E_{C,Z})=6$ and $h^0(H)=5$. Hence, the determinant map
$$\det:\bigwedge^2H^0(E_{C,Z})\longrightarrow H^0(H)$$
is not injective. Therefore, there exist two sections $s_1,s_2\in H^0(E_{C,Z})$ which are linearly independent, and $\det(s_1\wedge s_2)=0$. $s_1$ and $s_2$ form a sub-line bundle $L$ of $E_{C,Z}$ with $h^0(L)\geq2$. By taking the saturation of $L$, we can assume that $L$ is saturated. Since $(E_{C,Z}/L)^{\vee\vee}\cong H\otimes L^{\vee}$, we have $L.(H\otimes L^{\vee})\leq 3$ and, by Lemma 4.1, we have $L\in\mathcal{A}(H)$. Since $\Cliff(C)=1$, $L.(H\otimes L^{\vee})\geq3$ and hence, $L\in\mathcal{A}^0(H)$. Therefore, we have the following exact sequence.
$$0\longrightarrow L\longrightarrow E_{C,Z}\longrightarrow H\otimes L^{\vee}\longrightarrow 0$$
By Theorem 6.1 (i), $|H\otimes L^{\vee}|$ is an elliptic pencil and $\mathcal{O}_C(Z)=H\otimes L^{\vee}|_C$. Since, by Proposition 5.1 (ii), $L=F$ or $H\otimes F^{\vee}$, we have the assertion.$\hfill\square$

$\;$

\noindent {\bf{Acknowledgements}}

\smallskip

\smallskip

\noindent The author is partially supported by Grant-in-Aid for Scientific Research (16K05101), Japan Society for the Promotion of Science.

\end{document}